\documentclass[12pt]{amsart}
\usepackage{amsfonts,amscd,latexsym}

\makeatletter
\def\pf{%
  \par\topsep6\p@\@plus6\p@
  \trivlist
  \item[\hskip\labelsep\it\proofname.]\ignorespaces}

\makeatother

\numberwithin{equation}{section}

\newtheorem{thm}{Theorem}[section]
\newtheorem{lem}{Lemma}[section]
\newtheorem{cor}{Corollary}[section]

\theoremstyle{definition}

\theoremstyle{remark}
\newtheorem{rem}{Remark}

\begin{document}

\newcommand{\bfR}{\mathbb{ R}}
\newcommand{\bfC}{\mathbb{ C}}
\newcommand{\bfZ}{\mathbb{ Z}}
\newcommand{\bfH}{\mathbb{ H}}
\newcommand{\bfQ}{\mathbb{ Q}}
\newcommand{\bfN}{\mathbb{ N}}
\newcommand{\bfP}{\mathbb{ P}}
\newcommand{\bfT}{\mathbb{ T}}
\newcommand{\A}{\mathcal A}
\newcommand{\K}{\mathcal K}
\newcommand{\M}{\mathcal M}
\newcommand{\B}{\mathcal B}
\newcommand{\G}{\mathcal G}

\newcommand{\n}{\noindent}
\newcommand{\ds}{\displaystyle}
\newcommand{\e}{\epsilon}
\newbox\qedbox
\setbox\qedbox=\hbox{$\Box$}

\setlength{\baselineskip}{1.2\baselineskip}

\title[Isometric immersions]{Isometric immersions of space forms and soliton
theory}
\author[D. Ferus]{Dirk Ferus}
\address{Fachbereich Mathematik, MA 8-3 \\
Technische Universit\"at Berlin \\
Strasse des 17. Juni 136\\ 10622 Berlin, Germany}
\email[Dirk Ferus]{ferus@@sfb288.math.tu-berlin.de}
\thanks{Second author was partially supported by NSF grant DMS-9205293,
the SFB 288 at Technische Universit\"at Berlin and the
Graduierten Kolleg at Humboldt Universit\"at Berlin}
\author[F. Pedit]{Franz Pedit}
\address{Department of Mathematics \\
University of Massachusetts \\
Amherst, MA 01003}
\email[Franz Pedit]{franz@@gang.umass.edu}

\date{}

\maketitle

\section{Introduction}
The study of isometric immersions of space forms into space forms
is a classical problem of differential geometry. In its simplest form
it arises as the study of surfaces in $3$-space of constant (non-zero) Gaussian
curvature. In this case the integrability condition reduces to the sin resp.
sinh-Gordon equations. Due to the complicated structure of these equations
one focused over the past decades mainly on non-existence results
rather then the construction of explicit examples, the most famous one
being a result by Hilbert \cite{Hil} that there are no complete surfaces
of negative curvature in Euclidean $3$-space. A slightly different
spirit was present in the work of Bianchi and Darboux  (see also \cite{Ams})
where one finds many examples, some of which probably
have been forgotten over time and only recently been rediscovered in the
framework of what might be called {\em soliton geometry}. Having its
origin in the explicit description of all tori of constant mean curvature
\cite{PinSte,Boball} the finite gap integration scheme can also be
applied to obtain explicit parametrizations of a large class of
surfaces of constant Gaussian curvature \cite{MelSteLeeds, MelSte,
BobPindis, BobLeeds}.

In higher dimensions the situation is quite similar: there are a number
of interesting non-existence results \cite{Spi5, CheKui,Moo, Moopos,Ped,Gro}
for
isometric immersions
$ f: M^{m}(c)
\rightarrow \widetilde{M}^{n}(\widetilde{c}) $ between space forms
$M$ of dimension $m$, curvature $ c $, and space forms $ \widetilde{M}
^{n}(\widetilde{c}) $ of dimension $ n $, curvature $ \widetilde{c} $, but
no general scheme to construct such immersions. The integrability conditions
for
such an  immersion, the Gauss, Codazzi and Ricci equations, are sometimes
referred to as a generalized sin-Gordon system. Certain aspects
of its inverse scattering theory \cite{McK, Ter} as well as solutions obtained
by B\"acklund transformations (which generate a finite
dimensional solution space starting from a known, {\em trivial\/},
solution) \cite{TerTen} have recently been studied.

On the other hand, exterior differential calculus and the Cartan-K\"ahler
theory  \cite{GriJen} leads to a description of the space of local {\em real
analytic}
isometric immersions of space forms (in specific codimensions) in terms of
finitely
many functions of a single variable. In particular this
shows that this space is infinite dimesional.

In this paper we will develop the finite gap integration theory for
the integrability conditions for isometric immersions of space forms.
By this we mean that we will construct a hierarchy of commuting {\em finite
dimensional} ODE in Lax form whose solutions will give local real analytic
isometric
immersions of space forms. We call those isometric
immersions of {\em finite type}. The flows evolve on a certain
loop Lie algebra and can, at least in principle, be integrated by theta
functions
on the Jacobian of an algebraic curve, the {\em spectral curve} of the
Lax flows. In this sense the isometric immersion
equations for space forms are a completely integrable system and thus
can be regarded as a higher  dimesional (w.r.t. the flow variables) soliton
equation.
The explicit algebraic integration so far has only be
carried out in detail in the case of surfaces in $3$-space \cite{MelSte,
BobPindis} and the higher dimensional case would require substantial more work
on the
description of Riemann surfaces with two (commuting) involutions and certain
reality
conditions which will be done elsewhere. The initial condition for the Lax
flows is
a   matrix valued function of one variable so that our solution space
has at least functionally the correct dimension given by Cartan-K\"ahler theory
and we expect that an appropiate closure of the finite type solutions will
in fact account for all the local real analytic isometric immersions of space
forms.

In order to avoid degenerate cases we will make the following
assumptions on our isometric immersions $ f: M^{m}(c)\rightarrow \widetilde{M}
(\widetilde{c}) $ throughout the paper :
\renewcommand{\theenumi}{\roman{enumi}}
\begin{enumerate}
\item
$ c\neq 0 \neq \widetilde{c}, \quad c\neq \widetilde{c} $, \, and
\item
the normal bundle of $ f: M^{m}(c)\rightarrow \widetilde{M}
(\widetilde{c}) $ is flat.
\end{enumerate}
Note that for $ c < \widetilde{c} $ and $ n=2m-1 $ the second
assumption is always satisfied \cite{Moo} and, in this case, there are no
(local) isometric immersions for $ n\le 2m-2 $ \cite{Spi5}. It  is a classical
result due to Hilbert \cite{Hil, Spi4} that there are no complete isometric
immersions $ M^{2}(c)\rightarrow \widetilde{M}^{3}(\widetilde{c}) $
if $ c<\widetilde{c} $, $c<0$ and it is anybodies guess that this result
extends to complete isometric immersions $ M^{m}(c)\rightarrow
\widetilde{M}^{2m-1}(\widetilde{c}) $ for $ c<\widetilde{c} $, $c<0$.
Note that  for $c=0$ one always has the Clifford tori and $c>0$ cannot occur
due to the fact that such immersions induce global Chebycheff coordinates
\cite{Moo, Ped}.
In contrast, in the case $ c > \widetilde{c} $ one always has
the totally umbilical hypersurfaces. If $n\leq 2m-1$ and the immersion has
no umbilic points \cite{Moopos} then $n=2m-1$ and the normal bundle is again
flat.
Generally, flatness of the normal bundle is
a necessary assumption in higher codimensions to guarantee that
those isometric immersions are not too flabby to arise from completely
integrable systems \cite{FerPin}.  The case $ c=\widetilde{c} $
is treated exhaustively in \cite{FerPedflat}.

First we discuss how the structural equations for an isometric immersion
between
space forms can be rewritten as a  zero curvature condition involving an
auxiliary
({\em spectral\/}) parameter.  Thus each such isometric immersion is part of a
1-parameter family of isometric immersions and the deformation parameter turns
out to
be related to the induced curvature. We then express this fact as the flatness
condition on a loop algebra valued $1$-form with an algebraic constraint.
Such a reformulation of the equations at hand is
reminiscent in the theory of integrable systems and will be the starting point
for
their integration. To carry this through we  had to slightly modify
AKS or R-matrix theory \cite{BurPedLeeds}, since our equations could not be
treated
inside the standard setup. This modification was geometrically rooted since the
Lax
flows so obtained yield framings of the isometric
immersions which parallelize the normal bundle and thus are adapted to the
geometry of the situation.

A major part of this research was done while the second author was a
visiting member at the SFB 288 at Technische Universit\"at Berlin and enjoyed
the hospitality,
interest and encouragment present at this institution. Finally we would like
to thank Josef Dorfmeister for inspiring us to think about a modification
of standard AKS and R-matrix formalism.

\section{Associated family of isometric immersions}

Since our discussion will be local, we may (after scaling) assume
$ \widetilde{M}
^{n}(\widetilde{c}) $ to be either the Euclidean sphere $ S^{n}
\subset {\bfR}^{n+1} $ of curvature $ \widetilde{c}=1 $ or
hyperbolic space $ H^{n}\subset {\bfR}^{n+1}_{1} $ of curvature
$ \widetilde{c}=-1 $ realized as one sheet of the hyperboloid in
Lorentz space. We denote by $S$ either of these two standard
spaces.

Let $ f: M^{m}(c)\rightarrow S $ be an isometric immersion and let
$ F: M^{m}(c)\rightarrow \text{\bf SO}_{\epsilon}
(n+1) $ be an adapted framing, that is to say $ F_{0}=f, F_{1}, \dots, F_m $
are
tangential and $ F_{m+1}, \dots, F_n $ are normal. We denote by
\[A=F^{-1}dF\]
the pull back of the left Maurer-Cartan form of $ \text{\bf SO}_{\epsilon}
(n+1) $ by $ F $. (Here and in the sequel $\epsilon= \widetilde{c}=\pm
1$ and indicates whether one deals with the Euclidean or Lorentzian version of
the corresponding object).  Then
$ A $ is an $ \text{\bf so}_{\epsilon}(n+1)$-valued
1-form and since $ F $ is adapted,
\begin{equation} \label{eq:2.1}
A=\begin{pmatrix}
0 & -\e\theta^{T} & 0 \\
\theta & \omega & \beta \\
0 & -\beta^{T} & \eta
\end{pmatrix} ,
\end{equation}
where $ \theta $ is an $ {\bfR}^{m}$-valued 1-form on $ M^{m}
(c) $ (the dual frame to $ F_{1}, \dots, F_{m}) $, $ \omega $
is an $\text{\bf so}(m)$-valued 1-form (the Levi-Civita connection on
$ M^{m}(c) $), $ \eta $ is an $\text{\bf so}(n-m)$-valued 1-form (the
connection in the normal bundle) and $ \beta $ is the $ 2^{\text{nd}} $
fundamental form for $ f $.  The integrability conditions for the
existence of such a framing $ F $, the Maurer-Cartan equations
\begin{equation} \label{eq:2.2}
dA+A \wedge A=0\,,
\end{equation}
unravel to the fundamental equations for the isometric immersion
$ f $, the structure equation and the Gauss, Codazzi and Ricci
equations

\begin{subequations}{a} \label{eq:2.3}
\begin{gather}
d\theta+\omega \wedge \theta=0\,,\label{eq:2.3a}\\
d\omega + \omega \wedge \omega -\e\theta \wedge \theta^{T}-\beta
\wedge \beta^{T}=0\,,\label{eq:2.3b}\\
d\beta + \omega \wedge \beta + \beta \wedge \eta = 0\,,\label{eq:2.3c}\\
d\eta + \eta \wedge \eta-\beta^{T}\wedge \beta = 0\,.\label{eq:2.3d}
\end{gather}
\end{subequations}

In addition one has that the induced metric has curvature
$ c $,
\begin{equation} \label{eq:2.4}
d\omega + \omega \wedge \omega=c\,\theta \wedge \theta^T\,,
\end{equation}
and that the normal bundle is flat,
\begin{equation} \label{eq:2.5}
d\eta+\eta\wedge\eta=0\,.
\end{equation}
The starting point in soliton theory for integrating
the equations at hand is their reformulation as a zero curvature
condition (Maurer-Cartan equation) involving some auxiliary
(spectral) parameter.  Let $ A $ be an $ \text{\bf so}_{\e}(n+1)$-valued
1-form on $ M^{m}(c) $ of the form (\ref{eq:2.1}).  We define a
family of $ \text{\bf so}_{\e}(n+1, \bfC )$-valued 1-forms parametrized
by $ \lambda \in {\bfC}^* $ by
\begin{equation} \label{eq:2.6}
\widetilde{A}^{\lambda}=
\begin{pmatrix}
0 & \frac{-\e\sqrt{\e c}}{2}(\lambda + \lambda^{-1})\theta^T &0 \\
\frac{\sqrt{\e c}}{2}(\lambda + \lambda^{-1})\theta & \omega &
\frac{\sqrt{\e c}}{2\sqrt{1-\e c}}(\lambda-\lambda^{-1})\beta \\
0 & \frac{-\sqrt{\e c}}{2\sqrt{1-\e c}}(\lambda-\lambda^{-1})\beta^T
& \eta
\end{pmatrix}
\end{equation}
\begin{lem} \label{lem:Maurer-Cartan}
$ A $ solves (\ref{eq:2.3}), (\ref{eq:2.4}), and (\ref{eq:2.5})
if and only if $ \widetilde{A}^\lambda $ solves the Maurer-Cartan
equation
\[d\widetilde{A}^{\lambda}+\widetilde{A}^{\lambda}\wedge
\widetilde{A}^{\lambda}=0\]
for all $ \lambda $ (in a set with accumulation point) in ${\bfC}^{*}$.
\end{lem}

\begin{pf}
Since (\ref{eq:2.3a}) to (\ref{eq:2.3d}) are equivalent to the
Maurer-Cartan equation
(\ref{eq:2.2}) we have to verify that $ \widetilde{A}^\lambda $ solves
(\ref{eq:2.3a}) to (\ref{eq:2.3d}).  Since (\ref{eq:2.3a}) and
(\ref{eq:2.3c}) are homogeneous in
$ \theta $ and $ \beta $ (and neither $ \omega $ nor $ \eta $
involve $ \lambda $) they are trivially satisfied by $ \widetilde
{A}^\lambda $.  Also (\ref{eq:2.3d}) holds because of (\ref{eq:2.5}).
Now (\ref{eq:2.3b}) and (\ref{eq:2.4}) imply
\[(c-\epsilon) \theta\wedge\theta^T=\beta\wedge\beta^T=0\]
so that, using (\ref{eq:2.4}), we have
\[c\theta\wedge\theta^T-\frac{c}{2}(\lambda+\lambda^{-1})^2\theta
\wedge\theta^T-\frac{\epsilon c}{4(1-\epsilon c)}(\lambda-\lambda
^{-1})^2\beta\wedge\beta^T=0\,,\]
which shows that also (\ref{eq:2.3b}) holds for $ \widetilde{A}^\lambda $.
To prove the converse we compute the entries in $ d\widetilde{A}
^\lambda+\widetilde{A}^\lambda\wedge\widetilde{A}^\lambda=0 $ and
compare coefficients of equal powers of $ \lambda $:
\begin{align*}
&d(\frac{\sqrt{\epsilon c}}{2}(\lambda+\lambda^{-1})\theta)+
\omega\wedge(\frac{\sqrt{\epsilon c}}{2} (\lambda+\lambda^{-1})
\theta)=0\\
&d\omega+\omega\wedge\omega-\frac{c}{4}(\lambda+\lambda^{-1})^2\theta
\wedge\theta^T-\frac{\epsilon c}{4(1-\epsilon c)}(\lambda-\lambda^{-1})
^2\beta\wedge\beta^T=0\\
&d(\frac{\sqrt{\epsilon c}}{2\sqrt{1-\epsilon c}}(\lambda-\lambda
^{-1})\beta)+\omega\wedge(\frac{\sqrt{\epsilon c}}{2\sqrt{1-\epsilon
c}}(\lambda-\lambda^{-1})\beta)+(\frac{\sqrt{\epsilon c}}{2\sqrt
{1-\epsilon c}}(\lambda-\lambda^{-1})\beta)\wedge \eta=0\\
&d\eta+\eta\wedge\eta-\frac{c}{4(1-\epsilon c)}(\lambda-\lambda^{-1})
^2\beta^T\wedge\beta=0\,.
\end{align*}
The first and third identity yield (\ref{eq:2.3a}) and (\ref{eq:2.3c})
and the last identity gives (\ref{eq:2.5}) and (\ref{eq:2.3d}).
Expanding the second identity we obtain
\[d\omega+\omega\wedge\omega-\frac{c}{2}\theta\wedge\theta^T+\frac
{\epsilon c}{2(1-\epsilon c)}\beta\wedge\beta^T-(\lambda^2+\lambda
^{-2})(\frac{c}{4}\theta\wedge\theta^T+\frac{\epsilon c}{4(1-\epsilon
c)}\beta\wedge\beta^T)=0\]
and thus
\begin{align*}
&(1-\epsilon c)\theta\wedge\theta^T+\epsilon\beta\wedge\beta^T=0\\
&d\omega+\omega\wedge\omega=c\theta\wedge\theta^T
\end{align*}
which implies (\ref{eq:2.4}) and (\ref{eq:2.3b}).
\end{pf}
In the situation of Lemma~\ref{lem:Maurer-Cartan} we can integrate
$ (F^\lambda)^{-1}dF^\lambda=\widetilde{A}^\lambda $ to a (complex)
framing
\[F^\lambda: M^m\to \text{\bf SO}_\epsilon(n+1,\bfC)\]
for each $ \lambda $.  To obtain a real valued framing $ \widetilde{A}
^\lambda $ has to take values in $ \text{\bf so}_\epsilon(n+1) $.  This
is the case if and only if
\begin{equation}\label{eq:real}
\begin{tabular}{|c|c|c|} \hline
&
$ S=S^n $ &
$ S=H^n $ \\ \hline
$ \lambda $ & $ c^\lambda $ & $ c^\lambda $ \\ \hline
real & $ (0,1) $ & $ (-1,0) $ \\
imaginary & $ (-\infty,0) $ & $ (0,\infty) $ \\
unitary & $ (1,\infty) $ & $ (-\infty,-1) $ \\ \hline
\end{tabular}
\end{equation}

In either of these cases $ F^\lambda: M^m\to
\text{\bf SO}_\epsilon(n+1) $ is real and thus we obtain a family of
isometric immersions with
flat normal bundles (since (\ref{eq:2.5}) holds for each $ \widetilde
{A}^\lambda $)
\begin{equation}\label{eq:assfam}
f^\lambda=F^\lambda_0: M^m\to S
\end{equation}
from the first column of $ F^\lambda $.  The curvature $ c^\lambda $
of the induced metric is given by (\ref{eq:2.4}) expressed in the
coframe for $ f^\lambda $:
\[d\omega+\omega\wedge\omega=\epsilon\frac{4}{(\lambda+\lambda^{-1})^2}
\frac{\sqrt{c}}{2}(\lambda+\lambda^{-1})\theta\wedge\frac{\sqrt{c}}{2}
(\lambda+\lambda^{-1})\theta^T\,,\]
and thus
\begin{equation}\label{eq:asscurv}
c^{\lambda}=\epsilon\frac{4}{(\lambda+\lambda^{-1})^2}\,.
\end{equation}
Depending on the domain of $ \lambda $ the induced curvature
$ c^\lambda $ ranges over the intervals given in (\ref{eq:real}).
Finally, for
\begin{equation}\label{eq:lambda}
\lambda_0=\frac{1}{\sqrt{\e c}}(1+\sqrt{1-\epsilon c})
\end{equation}
we have $ \widetilde{A}^{\lambda_0}=A, f^{\lambda_0}=f $ and
$ c^{\lambda_0}=c $,
and thus recover the original immersion.  We summarize the discussion
so far in the following

\begin{lem}\label{lem:isoimm}
Let $ f: M^m\to S $ be an isometric immersion with
flat normal bundle,
$ F: M^m(c)\to \text{\bf SO}_\epsilon
(n+1) $ an adapted framing with induced Maurer-Cartan form $ A=F^{-1}dF $
and let $ \widetilde{A}^\lambda $ be given by (\ref{eq:2.6}).  Then
$ d\widetilde{A}^\lambda + \widetilde{A}^\lambda\wedge
\widetilde{A}^\lambda=0 $ for all $ \lambda $.

If $ \widetilde{A}
^\lambda $ satisfies the reality conditions (\ref{eq:real}) and
$ \widetilde{F}^\lambda $ integrates $ \widetilde{A}^\lambda $ then
$ \widetilde{F}^\lambda: M^m\to \text{\bf SO}_\epsilon(n+1) $
is an adapted framing for the isometric immersion $ f^\lambda=F^\lambda
_0: M^m(c^\lambda)\to S $ with flat normal bundle
and induced curvature $ c^\lambda=\epsilon\frac{4}{(\lambda+\lambda
^{-1})^2} $.  The original immersion $ f $ is recovered at
$ \lambda=\frac{1}{\sqrt{\epsilon c}}(1+\sqrt{1-\epsilon c}) $.
\end{lem}

Thus isometric immersions of space forms with flat normal bundle
come naturally in 1-parameter familes which we call the {\em
associated family}.

{}From Lemma~\ref{lem:Maurer-Cartan} we also obtain the converse to
Lemma~\ref{lem:isoimm}:

\begin{lem}\label{lem:ala}
Let
\begin{equation}\label{eq:ala}
\widetilde{A}^\lambda=\begin{pmatrix}
0 & -\epsilon (\lambda+\lambda^{-1})\theta^T &  0 \\
(\lambda+\lambda^{-1})\theta & \omega & (\lambda-\lambda^{-1})
\beta \\
0 & -(\lambda-\lambda^{-1})\beta^T & \eta
\end{pmatrix}
\end{equation}
be a family of $ \text{\bf so}_\epsilon(n+1)$-valued 1-forms on $ M^m $,
where the forms $ \theta $ and $ \beta $ may be imaginary
(to fulfill the reality conditions (\ref{eq:real}) ), satisfying
the Maurer-Cartan equation.  If $ \theta^T=(\theta^1,\dots,
\theta^m) $ are linearly independent then $ F^\lambda:
M^m\to \text{\bf SO}_\epsilon(n+1) $ integrating $ (F^\lambda)
^{-1}dF^\lambda=A^\lambda $ is an adapted framing for the
isometric immersion $ f^\lambda=F^\lambda_0: M^m(c^\lambda)
\to S $ with induced metric $ c^\lambda=\epsilon
\frac{4}{(\lambda+\lambda^{-1})^2} $ and flat normal bundle.
\end{lem}

The gist of this reformulation is that the construction of
isometric immersions of space forms with flat normal bundle is
equivalent to the construction of a certain family of
$ \text{\bf so}_\epsilon(n+1) $ valued 1-forms (\ref{eq:ala}) satisfying
the Maurer-Cartan equation.  Notice that such a reformulation
is well known in the theory of harmonic maps of Riemann
surfaces into Lie groups and symmetric spaces \cite{BurFerPedPin}.

\section{Loop algebra formulation}

In the previous section we discussed how the equations for an
isometric immersion $ f: M^m\to S $ with flat
normal bundle can be written as the zero curvature condition
for a family (``loop'') of Lie algebra valued 1-forms
(\ref{eq:ala}).  The proper setting for this are loop algebras.
Let $ \mathfrak{g}=\text{\bf so}_\epsilon(n+1) $ with complexification
$ \mathfrak{g}^{\bfC}= \text{\bf so}_\epsilon(n+1,\bfC)$ and define the loop
algebra
\[
\Lambda\mathfrak{g}^{\bfC}=\{\xi:{\bfC}^{*}\to
\mathfrak{g}^{\bfC}\,; \xi \text{ polynomial in}\; \lambda\; \text{and}\;
\lambda^{-1}\}
\]
which is a (complex) Lie algebra under pointwise bracket.
If $ C $ denotes either of $ {\bfR}^*, i{\bfR}^* $ or $ S^1 $
we have the real subalgebras (c.f. also (\ref{eq:real}))
\[
\Lambda\mathfrak{g}=\{\xi: C\to \mathfrak{g}\}
\subset \Lambda \mathfrak{g}^{\bfC}
\]
corresponding to the conjugations
\begin{equation}\label{eq:reality}
\bar{\xi}(\lambda)=\overline{\xi(\bar{\lambda})}\,,\qquad
\bar{\xi}(\lambda)=\overline{\xi(-\bar{\lambda})}\,,\qquad
\bar{\xi}(\lambda)=\overline{\xi(1/\bar{\lambda})}\,.
\end{equation}
Note that
\[\widetilde{A}^\lambda =\, \begin{pmatrix}
0 & -\epsilon(\lambda+\lambda^{-1})\theta^T & 0 \\
(\lambda+\lambda^{-1})\theta & \omega & (\lambda-\lambda^{-1})
\beta \\
0 & -(\lambda-\lambda^{-1})\beta^T & \eta
\end{pmatrix}\]

for $\lambda\in C$ can be regarded as a $ \Lambda\mathfrak{g}$-valued 1-form
with
the
following symmetries:
\begin{align*}
& \widetilde{A}^{-\lambda}=Ad\,P \widetilde{A}^\lambda \\
& \widetilde{A}^{1/\lambda} = Ad\,Q \widetilde{A}^\lambda
\end{align*}
where
\[P=\biggl(\begin{tabular}{c|c|c}
-1 &  \\ \hline
& 1 \\ \hline
& & -1
\end{tabular}\biggr)
\text{ and }
Q=\biggl(\begin{tabular}{c|c|c}
1 & \\ \hline
& 1 \\ \hline
& & -1
\end{tabular}\biggr).\]
This motivates to consider the following involutions on
$ \Lambda \mathfrak{g} $:
\begin{align*}
& (\sigma\xi)(\lambda)=Ad\,P\xi(-\lambda) \\
& (\tau\xi)(\lambda)=Ad\,Q\xi(1/\lambda) \,.
\end{align*}
We let
\[\Lambda\mathfrak{g}_{\sigma,\tau}=\{\xi\in\Lambda\mathfrak{g};
\sigma\xi=\xi, \tau\xi=\xi\}\]
be the subalgebra fixed under $ \sigma $ and $ \tau $.

\begin{lem}\label{lem:frakg}
Let  $ \xi=\sum_{k\in{\bfZ}}\lambda^k\xi_k\in\Lambda\mathfrak{g} $.
Then $\xi\in\Lambda\mathfrak{g}_{\sigma,
\tau} $ if and only if $ \xi_0\in\mathfrak{g}^P\cap\mathfrak{g}^Q,
\xi_{\text{even}}\in\mathfrak{g}^P, \xi_{\text{odd}}\in\mathfrak{g}
^{-P} $ and $ \xi_{-k}=Ad\,Q\xi_k $.  (Here $ \mathfrak{g}^{\pm
P}, \mathfrak{g}^{\pm Q} $ denote the $ \pm 1$-eigenspaces of
$ Ad\,P $ resp. $ Ad\,Q$ on $\mathfrak{g} $).
\end{lem}

\begin{pf}
This follows at once by comparing coefficients in
$ \sigma\xi=\tau\xi=\xi $.
\end{pf}

For $d\in \bfN$ let
\[
\Lambda_d=\{\xi\in\Lambda\mathfrak{g}_{\sigma,\tau};
\xi=\sum_{|k|\leq d}\lambda^k\xi_k\}
\]
be the subspace of Laurent polynomial loops of degree at most $d$.

\begin{cor}\label{cor:lambda}
$ \xi\in\Lambda_1 $ if and only if
\[\xi=\biggl(\begin{array}{c|c|c}
0 & -\epsilon (\lambda+\lambda^{-1})a^T & 0 \\ \hline
(\lambda+\lambda^{-1})a & A & (\lambda-\lambda^{-1})C \\ \hline
0 & -(\lambda-\lambda^{-1})C^T & B
\end{array}\biggr) \]
\end{cor}

\begin{pf}
This follows from Lemma~\ref{lem:frakg} together with the fact that

\begin{equation}\label{eq:fixp}
\mathfrak{g}^P=\biggl\{\biggl(\begin{tabular}{c|c|c}
{*}& & * \\ \hline
& * &\\ \hline
{*} & & *
\end{tabular}\biggr)\biggr\},
\quad
\mathfrak{g}^{-P}=\biggl\{\biggl(\begin{tabular}{c|c|c}
& * &  \\ \hline
{*} & & * \\\hline
& * &
\end{tabular}\biggr)\biggr\}
\end{equation}
\begin{equation}\label{eq:fixq}
\mathfrak{g}^Q=\biggl\{\biggl(\begin{tabular}{c|c|c}
{*}& * & \\ \hline
{*} & * & \\ \hline
& & *
\end{tabular}\biggr)\biggr\},
\quad
\mathfrak{g}^{-Q}=\biggl\{\biggl(\begin{tabular}{c|c|c}
& & * \\ \hline
& & * \\ \hline
* & * &
\end{tabular}\biggr)\biggr\}
\end{equation}
\end{pf}

Corollary~\ref{cor:lambda}
together with the reformulation of the isometric immersion
equations in \linebreak[4]
Lemma~\ref{lem:isoimm} and Lemma~\ref{lem:ala} now
yield

\begin{cor}\label{cor:natural}
There is a natural correspondence between isometric immersions
$ f: M^m\to S $ of space forms with flat
normal bundle and $ \Lambda_1$-valued 1-forms $ \widetilde{A}
: TM\to \Lambda_1 $ satisfying the
Maurer-Cartan equation (and whose first row
$ \widetilde{A}_{0,-}=(0,\epsilon(\lambda+\lambda^{-1})\theta^T,0) $
has $ \theta^T=(\theta^1, \dots ,\theta^m) $ linearly independent).
\end{cor}

The construction of flat loop algebra valued 1-forms satisfying
an algebraic constraint (i.e. taking values in $ \Lambda_1 $)
is reminiscent in the theory of soliton equations: they appear
as solutions to certain completely integrable Lax-type equations
\cite{BurPedLeeds}.

\section{Integration of the isometric immersion equations for
space forms}

We are now going to discuss a recipe for the construction of
$ \Lambda_1$-valued flat 1-forms on $ {\bfR}^m $ from finite
dimensional commuting Lax flows.  Following existing
nomenclature we will call the isometric immersions so obtained
{\em finite type isometric immersions}.  Such immersions will
be real analytic by construction (in fact given by theta
functions) and thus cannot account for all isometric immersions
of space forms.  Our solutions will be parametrized by a finite
number of functions in one variable which relates to well known
results obtained by Cartan-K\"ahler-theory \cite{GriJen}.  We begin by
recalling
some facts about the integration of certain Lax equations on Lie
algebras. Let $\G$ be a Lie algebra (in our case a loop algebra)
which has a vector space direct sum decomposition
\[{\mathcal G}={\mathcal P}\oplus {\mathcal A} \oplus {\mathcal M} \]
with
\[{\mathcal K} = {\mathcal A} \oplus {\mathcal P},\qquad {\mathcal B} =
{\mathcal A}
 \oplus {\mathcal M}
\qquad \text{and}\;\;{\mathcal A}\]
Lie subalgebras and commutation relations
\[
[{\mathcal A},{\mathcal P}]\subseteq{\mathcal P}, \qquad [{\mathcal
A},{\mathcal
M}]\subseteq{\mathcal M}\,.
\]
Denote the corresponding projections by $\pi_{\mathcal P}$,
$\pi_{\mathcal A}$ and
$\pi_{\mathcal M}$.  An ad-equivariant vector field on $\mathcal G$ is a map
\[V: \mathcal G \rightarrow \mathcal G\]
satisfying
\begin{equation}\label{eq:adequ}
d_\xi V[\xi,\eta]=[V(\xi),\eta]
\end{equation}
for $\xi,\eta\in\mathcal G$. Natural examples
of such maps are the gradients (with respect to an invariant inner product) of
ad-invariant functions on $\mathcal G$. Notice that (\ref{eq:adequ}) implies
that two ad-equivariant vector fields $V,{\widetilde V}:{\mathcal
G}\to{\mathcal G}$
commute pointwisely, i.e., that
\begin{equation}\label{eq:com}
[V(\xi),{\widetilde V}(\xi)]=0
\end{equation}
for $\xi\in{\mathcal G}$. Given an
ad-equivariant vector field
$V$ on $\mathcal G$ we define a vector field $X^V$ on $\mathcal K$ by
\[X^V(\xi)=[\xi, \pi_{\mathcal P}V(\xi)], \quad \xi\in\mathcal K\,.\]
In order to make sense of the derivative of a map on
$\G$ without specifying a topology we assume that the images under the map of
finite dimensional (vector) subspaces are contained in finite dimensional
(vector)
subspaces. Note that if $V$ has this property so automatically has $X^{V}$.
\begin{lem}\label{lem:com}
Let $V,{\widetilde V}:{\mathcal G}\to{\mathcal G}$ be ad-equivariant vector
fields.
Then, for $\xi\in{\mathcal K}$,
\[
[X^{V}, X^{{\widetilde V}}]_{C^{\infty}}(\xi)=[\xi,\pi_{\mathcal
A}[\pi_{\mathcal
P}V(\xi),\pi_{\mathcal P}{\widetilde V}(\xi)]]\,.
\]
\end{lem}
The proof consists of a straightforward calculation using (\ref{eq:adequ}) and
(\ref{eq:com}). Observe that in the case ${\mathcal A}=0$, i.e., when
${\mathcal G}=
{\mathcal K}\oplus{\mathcal B}$ is the (vector space) direct sum of two Lie
subalgebras, the vector fields $X^{V}$ and $X^{{\widetilde V}}$ commute, and
one is
in the  realm of standard AKS resp. R-matrix theory \cite{AdlMor, BurPedLeeds}.

\begin{lem}\label{lem:flat}
Let $V_1,\dots,V_m$ be ad-equivariant vector fields on $\mathcal G$ and assume
that
for all $\xi\in{\mathcal K}$
\begin{equation}\label{eq:cond}
\pi_{\mathcal A}[\pi_{\mathcal P}V_{i}(\xi),\pi_{\mathcal P}V_{j}(\xi)]=0\,.
\end{equation}
Then the system of ODE
\[
d\xi=\sum^m_{i=1}X^{V_{i}}(\xi)dx^i=\sum^m_{i=1}[\xi, \pi_{\mathcal
P}V_{i}(\xi)]dx^i\,,
\quad
\xi(0)\in \mathcal K
\]
has a unique (local) solution $\xi: U\subset\bfR^m\rightarrow \mathcal K$, and
the
$\mathcal P$-valued $1$-form on $U$,
\[
\widetilde{A}=\sum^m_{i=1}\pi_{\mathcal P}V_i(\xi)dx^i\,,
\]
satisfies the Mauer-Cartan equation
\[
d\widetilde{A}+\frac{1}{2}[\widetilde{A}\wedge\widetilde{A}]=0\,.
\]
\end{lem}
\begin{pf}
By Lemma~\ref{lem:com} the vector fields $X^{V_{i}}$ on $\mathcal K$ commute.
Thus the system
\[
d\xi=\sum^m_{i=1}X^{V_{i}}(\xi)dx^i
\]
is well-defined and so has a unique solution to any initial condition
$\xi(0)\in{\mathcal K}$. The final statement follows from an analogous
calculation
as in the proof of Lemma~\ref{lem:com} which gives
\[
d\widetilde{A}+\frac{1}{2}[\widetilde{A}\wedge\widetilde{A}]=
\sum^m_{i,j=1}\pi_{\mathcal A}[\pi_{\mathcal P}V_{i}(\xi),\pi_{\mathcal
P}V_{j}(\xi)]dx^i\wedge
dx^j=0\,.
\]
\end{pf}
\begin{rem}
The previous two Lemmas are simple modifications of standard results in
AKS resp. R-matrix theory. Even though one  cannot expect  the commutativity
conditions (\ref{eq:cond}) to hold for general ad-equivariant vector fields
there
are interesting geometric applications where they do in fact hold. One of them
are the isometric immersion equations (Corollary~3.2), which do not fit into
the standard scheme, but can be treated with this more general setup. A more
detailed
study of the abstract general situation will perhaps be done elsewhere.
\end{rem}

To apply the above considerations to the integration of the isometric
immersion equations we first work in the complex setup and put
\begin{align*}
{\mathcal
G}&=\Lambda\mathfrak{g}^{\bfC}_{\sigma}=\{\xi\in\Lambda\mathfrak{g}^{\bfC};
\sigma\xi=\xi\} \\
\intertext{and}
{\mathcal
K}&=\Lambda\mathfrak{g}^{\bfC}_{\sigma,\tau}=
\{\xi\in\Lambda\mathfrak{g}^{\bfC};
\sigma\xi=\tau\xi=\xi\}\,.
\end{align*}
A natural vector space complement $\M$ to $\K$ in
$\Lambda\mathfrak{g}^{\bfC}_{\sigma}$
is given as follows: note that
\begin{equation}\label{eq:pm}
\G=\G_{-}\oplus\G_{0}\oplus\G_{+}
\end{equation}
with
\[
\G_{\pm}=
\{\xi\in\Lambda\mathfrak{g}^{\bfC}_{\sigma};\xi=
\sum_{k><0}\lambda^{k}\xi_{k}\}\,,
\qquad \G_{0}=({\mathfrak g}^{\bfC})^{P}\,.
\]
Given $\xi\in\Lambda\mathfrak{g}^{\bfC}_{\sigma}$ we write
\[
\xi=\xi_{-}+\xi_{0}+\xi_{+}
\]
according to (\ref{eq:pm}) and decompose
\begin{equation}\label{eq:km}
\xi=(\xi_{-}+\tau\xi_{-}+\frac{\xi_{0}+\tau\xi_{0}}{2})+(\xi_{+}-\tau\xi_{-}+
\frac{\xi_{0}-\tau\xi_{0}}{2})
\end{equation}
to obtain the factorization
\[
\G=\K\oplus\M
\]
where
\[
\M=({\mathfrak g}^{P}\cap{\mathfrak g}^{-Q})^{\bfC}\oplus\G_{+}\,.
\]
Unfortunately $\M$ is not a Lie subalgebra since
\[
0\neq [{\mathfrak g}^{P}\cap{\mathfrak g}^{-Q},{\mathfrak g}^{P}\cap{\mathfrak
g}^{-Q}]\subseteq{\mathfrak g}^{P}\cap{\mathfrak g}^{Q}\,.
\]
In fact,
\[
{\mathfrak g}^{P}\cap{\mathfrak g}^{Q}=\text{\bf
so}_\epsilon(1+m)\oplus\text{\bf
so}(n-m)=:{\mathfrak a}_{1}\oplus{\mathfrak a}_{2}=\biggl\{\biggl(
\begin{tabular}{c|c|c} & &  \\ \hline &* &\\ \hline & & *
\end{tabular}\biggr)\biggr\}
\]
and
\[
[{\mathfrak g}^{P}\cap{\mathfrak g}^{-Q},{\mathfrak g}^{P}\cap{\mathfrak
g}^{-Q}]\subseteq{\mathfrak
a}_{2}\,.
\]
Thus, setting
\begin{align*}
\A&={\mathfrak a}_{2}^{\bfC}=\biggl\{\biggl(\begin{tabular}{c|c|c} & &  \\
\hline &
&\\ \hline & & *
\end{tabular}\biggr)\biggr\}\,,\\
{\mathcal P}&=\{\xi\in\K\,;\,\xi_{0}\in{\mathfrak a}_{1}^{\bfC}\}\,,
\end{align*}
we arrive at
\[
\G={\mathcal P}\oplus\A\oplus\M
\]
with
\[
\K={\mathcal P}\oplus\A\,,\qquad\B=\A\oplus\M\qquad\text{and}\;\;\A
\]
Lie subalgebras satisfying the commutation relations
\[
[\A,{\mathcal P}]\subseteq{\mathcal P}\,,\qquad [\A,\M]\subseteq\M\,.
\]
Finally we introduce the ad-equivariant vector fields relevant to our
situation.

\begin{lem}\label{lem:setup}
(i) The map $V: \mathcal G\rightarrow \mathcal G$ defined by
\[V(\xi)=\lambda^{2k}\xi^{2\ell-1}, \qquad k\in{\bfZ}, \ell\in{\bfN}\]
is an ad-equivariant vector field on $\mathcal G$.

(ii) Let $d\in{\bfN}$ be odd and define
$V_{\ell}(\xi)=\lambda^{d(2\ell-1)-1}\xi^{2\ell-1}, \ell\in{\bfN}$.  Then the
corresponding vector fields on $\mathcal K$
\[
X^{V_{\ell}}(\xi)=[\xi, \pi_{\mathcal P}V_{\ell}(\xi)]\,,
\]
are tangential to the finite dimensionsal subspace
$\Lambda{_d}^{\bfC}=\{\xi\in\K\,;\,\xi=\sum_{|k|\leq
d}\lambda^{k}\xi_{k}\,\}\subset
\K$ and
\[
\pi_{\mathcal P}V_{\ell}(\xi)\in \Lambda_1^{\bfC}\,.
\]

(iii) $[X^{V_{i}},X^{V_{j}}]_{C^{\infty}}=0$, i.e., the vector fields
$X^{V_{\ell}}$
commute.

(iv) (Reality conditions) If $\xi\in\Lambda_{d}$ is real (w.r.t. any of
the three reality conditions
(\ref{eq:reality}) $\lambda\in\bfR^{*}$,$i\bfR^{*}$ or $S^1$) then also
$\pi_{\mathcal P}V_{\ell}(\xi)\in \Lambda_1$ is real (and thus the vector
fields
$X^{V_{\ell}}$ are real, i.e., remain tangent to '$\Lambda_d$).

\end{lem}

\begin{pf}
(i) Since $\xi\in{\mathcal G}=\Lambda\mathfrak{g}_{\sigma}^{\bfC}$ we have
\begin{align*}
\begin{split}
\sigma V(\xi)(\lambda)&=(-\lambda)^{2k}Ad P\xi(-\lambda)^{2\ell-1}=
\lambda^{2k}(Ad P\xi(-\lambda))^{2\ell-1} \\
&=\lambda^{2k}\xi(\lambda)^{2\ell-1}=V(\xi)(\lambda)
\end{split}
\end{align*}
so that $ V:\mathcal G\to\mathcal G$ is well-defined.  To verify
(\ref{eq:adequ})
we simply compute
\[d_\xi V[\xi,\eta]=\lambda^{2k}\sum_{i+j=2\ell-2}\xi^i[\xi,\eta]
\xi^j=\lambda^{2k}(\xi^{2\ell-1}\eta-\eta\xi^{2\ell-1})=[V(\xi),
\eta]\,.\]

(ii) Let $\xi=\sum_{|j|\leq d}\lambda^j\xi_j\in\Lambda_{d}^{\bfC} $ then
\[
V_{\ell}(\xi)=\lambda^{d(2\ell-1)-1}\xi^{2\ell-1}=\lambda^{-1}\xi
^{2\ell-1}_{-d}+\sum^{2d(2\ell-1)-1}_{j=0}\lambda^j\widetilde{\eta}_{j}\,,
\]
where
\begin{equation}\label{eq:const}
\widetilde{\eta}_{0}=\sum_{k=0}^{2l-2}(\xi_{-d})^{k}\xi_{-d+1}(\xi_{-d})^{2l-2-k}\,.
\end{equation}
Thus, by (\ref{eq:km}),
\begin{equation}\label{eq:lam1}
\pi_{\mathcal P}V_{\ell}(\xi)=\lambda^{-1}\xi^{2\ell-1}_{-d}+\eta_0+\lambda
AdQ\xi^{2\ell-1}_{-d}\in \Lambda_1^{\bfC}
\end{equation}
with
$\eta_0=\pi_{{\mathfrak a}_{1}^{\bfC}}
\frac{\widetilde{\eta}_{0}+AdQ\widetilde{\eta}_{0}}{2}\in{\mathfrak
a}_{1}^{\bfC}$.
{}From this it is clear that
\[
X^{V_{\ell}}(\xi)=[\xi,\pi_{\mathcal P}V_{\ell}(\xi)]\in \Lambda_{d}^{\bfC}\,,
\]
i.e., $ X^{V_{\ell}}$ is tangential to $\Lambda_{d}^{\bfC}$.

(iii) Due to Lemma~\ref{lem:com} we have to verify that for
$\xi\in\Lambda_{d}^{\bfC}$
\[
\pi_{\A}[\pi_{\mathcal P}V_{i}(\xi),\pi_{\mathcal P}V_{j}(\xi)]=0\,.
\]
{}From (\ref{eq:lam1}) and the fact that $\eta_0\in{\mathfrak a_{1}}^{\bfC}$ it
suffices to
show that
\[
\pi_{\A}[\xi_{-d}^{2i-1}, AdQ\xi_{-d}^{2j-1}]=0\,.
\]
But this last is seen from the more general fact that
\[
\pi_{\A}[X,Y]=0
\]
for $X,Y\in({\mathfrak g}^{-P})^{\bfC}$ implies
\[
\pi_{\A}[X, AdQ\,Y]=0\,,
\]
which follows at once from (\ref{eq:fixp}) and the specific form of $Q$.

(iv) If $\lambda$ is real then the statement is obvious. In case $\lambda$ is
purely imaginary $\xi=\sum_{|k|\leq d}\lambda^{k}\xi_{k}\in\Lambda_{d}$
is equivalent to $\overline{\xi_{k}}=(-1)^k\xi_{k}$. From this and
(\ref{eq:lam1}), (\ref{eq:const}) it is clear that $\pi_{\mathcal
P}V_{\ell}(\xi)$
is real. Finally, if $\lambda$ is unitary the reality condition
(\ref{eq:reality})
for $\xi\in\Lambda_{d}$ gives $\xi_{-k}=\overline{\xi_{k}}=AdQ\xi_k$ so that
by (\ref{eq:lam1}), (\ref{eq:const}) we again conclude that
$\pi_{\mathcal P}V_{\ell}(\xi)$ is real.
\end{pf}

Putting together the above discussion with the results
in the previous sections we obtain a recipe for the construction
of (local) isometric immersions of space forms into $S^n$ and
$H^n$ from a hierarchy of finite dimensional ODE.

\begin{thm}\label{thm:recipe}
Let $ d\in{\bfN}$ be odd, $m:=[\frac{n+1}{2}]$ and recall
the ad-equivariant vector fields $V_{\ell}(\xi)=\lambda^{d(2
\ell-1)-1}\xi^{2\ell-1}, \ell=1,\dots,m$.

(i) The system of ODE
\begin{equation}\label{eq:laxpair}
d\xi=[\xi,\sum^m_{\ell=1}\pi_{\mathcal P}V_{\ell}(\xi)dx^{\ell}],\quad \xi
(0)=\overset{\circ}{\xi}\in\Lambda_d
\end{equation}
has a unique (local) solution $ \xi: U\subset{\bfR}^m\to
\Lambda_d$ for any initial condition $\overset{\circ}{\xi}\in
\Lambda_d$.

(ii) If $ \xi: U\rightarrow \Lambda_d$ is a solution to
(\ref{eq:laxpair}) then
\[\widetilde{A}=\sum^m_{\ell=1}\pi_{\mathcal P}V_{\ell}(\xi)dx^\ell\]
is a $\Lambda_1$-valued $1$-form on $U\subset{\bfR}^m$ solving
the (matrix) Maurer-Cartan equation
\[d\widetilde{A}+\widetilde{A}\wedge\widetilde{A}=0\,.\]
Thus, integrating $(F^\lambda)^{-1}dF^\lambda=\widetilde{A}^\lambda,
F^\lambda(0)=1$, to a framing $F^{\lambda}:U\to \text{\bf SO}_{\epsilon}(n+1)$
gives a family of isometric immersion
\begin{equation}\label{eq:immers}
f^\lambda=F^\lambda_0: U\to S=S^n \text{
or } H^n
\end{equation}
with induced curvature $ c^\lambda=\epsilon\frac{4}{(\lambda+\lambda
^{-1})^2}$.
\end{thm}

\begin{pf}
(i) follows from Lemma~\ref{lem:flat} together with Lemma~\ref{lem:setup}. (ii)
is a consequence of the description of isometric immersions as special loops
of flat $\text{\bf so}_{\epsilon}(1+n)$-valued $1$-forms, c.f. Lemma~2.3 and
Corollary~3.2, together with Lemmas \ref{lem:setup} and \ref{lem:flat}.
\end{pf}

\begin{rem}
(i) Due to the fact that
\[
\widetilde{A}=\sum^m_{\ell=1}\pi_{\mathcal P}V_{\ell}(\xi)dx^\ell\in
\biggl\{\biggl(\begin{tabular}{c|c|c} &*&  \\ \hline {*}& * & *\\ \hline &*&
\end{tabular}\biggr)\biggr\}
\]
the framing constructed from the Lax flows (\ref{eq:laxpair}) parallelizes the
normal bundle of the isometric immersion (\ref{eq:immers}) which makes the
flows
geometrically adapted.

(ii) Since the ad-equivariant vector fields $V_{\ell}$ arise from (shifts of)
gradients of ad-invariant functions on the finite dimensional Lie algebra
$\text{\bf so}_{\epsilon}(1+n)$ only rank $= [\frac{n+1}{2}]$ many vector
fields are independent, i.e., one cannot construct in this way isometric
immersions
of space forms of dimension exceeding $m=[\frac{n+1}{2}]$. This fact is also
geometrically rooted: it is well-known that there are no isometric immersions
$f:M^{m}(c)\to\tilde{M}^{n}(\tilde{c})$ in case $c<\tilde{c}$ and $n\leq 2m-2$.

(iii) Using Cartan-K\"ahler theory one can show that local real analytic
isometric immersions
$f:M^{m}(c)\to\tilde{M}^{2m-1}(\tilde{c})$ with $c<\tilde{c}$ depend
on finitely many functions in one variable. Our scheme produces
such real analytic isometric immersions from an arbitrary initial
condition $\overset{\circ}{\xi}\in\Lambda_d$, $d\in\bfN$ odd, which
is a certain $\text{\bf so}_{\epsilon}(1+n)$-valued Laurent polynomial (of
arbitrary odd degree) in the variable $\lambda$.
This indicates that an appropriate closure of the solutions so constructed
should account for all real analytic solutions, but at present we have little
idea how to adress this issue.

(iv) If the target space form is the sphere $S^n$ and hence the corresponding
Lie
algebra
$\mathfrak g =\text{\bf so}_{\epsilon}(1+n)$ compact, then there exists an
ad-invariant
positive definite (appropriately weighted $L^2$) inner product on
$\Lambda\mathfrak g$
which makes the flows (\ref{eq:laxpair}) evolve on Euclidean spheres in
$\Lambda_d$ and thus complete. But in general we cannot show that the flows
(\ref{eq:laxpair}) have global, i.e., defined on all of $\bfR^{m}$, solutions,
in
which case there still remains the question of whether the isometric immersion
(\ref{eq:immers}) is globally defined or whether its differential
drops rank. This is of course related to the conjecture that
there are no complete isometric immersions
$f:M^{m}(c)\to\tilde{M}^{2m-1}(\tilde{c})$ for $c<\tilde{c}$, $c<0$
which is a classical Theorem first proved by Hilbert in the case
of surfaces in Euclidean $3$-space.
\end{rem}

\bibliographystyle{amsplain}

\ifx\undefined\bysame
\newcommand{\bysame}{\leavevmode\hbox to3em{\hrulefill}\,}
\fi

\end{document}